\input harvmac
\newcount\figno
\figno=0
\def\fig#1#2#3{
\par\begingroup\parindent=0pt\leftskip=1cm\rightskip=1cm\parindent=0pt
\global\advance\figno by 1
\midinsert
\epsfxsize=#3
\centerline{\epsfbox{#2}}
\vskip 12pt
{\bf Fig. \the\figno:} #1\par
\endinsert\endgroup\par
}
\def\figlabel#1{\xdef#1{\the\figno}}
\def\encadremath#1{\vbox{\hrule\hbox{\vrule\kern8pt\vbox{\kern8pt
\hbox{$\displaystyle #1$}\kern8pt}
\kern8pt\vrule}\hrule}}

\overfullrule=0pt

%
\def\underarrow#1{\vbox{\ialign{##\crcr$\hfil\displaystyle
 {#1}\hfil$\crcr\noalign{\kern1pt\nointerlineskip}$\longrightarrow$\crcr}}}
%


%

\def\inbar{\vrule height1.5ex width.4pt depth0pt}
\def\IC{\relax\hbox{\kern.25em$\inbar\kern-.3em{\rm C}$}}
\def\IR{\relax\hbox{\kern.25em$\inbar\kern-.3em{\rm R}$}}
\def\IZ{\relax\ifmmode\hbox{Z\kern-.4em Z}\else{Z\kern-.4em Z}\fi}

\font\zfont = cmss10 

\def\bigone{\hbox{1\kern -.23em {\rm l}}}
\def\ZZ{\hbox{\zfont Z\kern-.4emZ}}


\def\drawbox#1#2{\hrule height#2pt
        \hbox{\vrule width#2pt height#1pt \kern#1pt
              \vrule width#2pt}
              \hrule height#2pt}

\def\Asym#1#2{\vcenter{\vbox{\drawbox{#1}{#2}
              \kern-#2pt       
              \drawbox{#1}{#2}}}}

\batchmode
  \font\bbbfont=msbm10
\errorstopmode
\newif\ifamsf\amsftrue
\ifx\bbbfont\nullfont
  \amsffalse
\fi
\ifamsf
\def\IR{\hbox{\bbbfont R}}
\def\IC{\hbox{\bbbfont C}}

\def\IZ{\hbox{\bbbfont Z}}

\def\IQ{\hbox{\bbbfont Q}}
\def\IN{\hbox{\bbbfont N}}
\def\IK{\hbox{\bbbfont K}}

\def\qed{\vbox{\hrule\hbox{\vrule\kern3pt
\vbox{\kern6pt}\kern3pt\vrule}\hrule}}


\midinsert
\endinsert


\nref\heis{Heisenberg W 1927 {\it Z. Phys.} {\bf 43} 172}

\nref\ken{Kennard E H 1927 {\it Z. Phys.} {\bf 44} 326}

\nref\weyl{Weyl H 1928 {\it Gruppentheorie und Quantenmechanik}
(Leipzig: Hirzel Verlag) }

\nref\rob1{Robertson H P 1929 {\it Phys. Rev.} {\bf 34} 163}

\nref\rob2{Robertson H P 1930 {\it Phys. Rev.} {\bf 35} 667A}

\nref\schr{Schr\"odinger E 1930 {\it Sitz. Preus. Acad. Wiss.
(Phys.-Math. Klasse)} 296}

\nref\rob3{Robertson H P 1934 {\it Phys. Rev.} {\bf 46} 794}

\nref\dodon{Dodonov V V, Kurmyshev E V and Man'ko V I {1980} {\it Phys.
Lett.} {\bf 79A} 150}

\nref\trif1{Trifonov D A 1993 {\it J. Math. Phys.} {\bf 34} 100}

\nref\trif1a{Trifonov D A 1994 {\it J. Math. Phys.} {\bf 35} 2297}

\nref\Sudar{Sudarshan E C G, Chiu C B and Bhamathi G 1995 {\it Phys. Rev. 
A} {\bf 52} 43}

\nref\trif2{Trifonov D A 1997 {\it J. Phys. A: Math. Gen.} {\bf 30} 5941}

\nref\trif3{Trifonov D A, Donev S G 1998 {\it J. Phys. A: Math. Gen.} {\bf
31} 8041}

\nref\trif4{Trifonov D A 2000 {\it J. Phys. A: Math. Gen.} {\bf
33} L299}

\nref\trif5{Trifonov D A 2001 {\it J. Phys. A: Math. Gen.} {\bf
34} L75}

\nref\trif6{Trifonov D A 1999 The uncertainty way of generalization of
coherent states {\it Preprint} quant-ph/9912084 v5}

\nref\trif7{Trifonov D A 2001 Generalizations of Heisenberg uncertainty
relation {\it Preprint} quant-ph/0112028}

\nref\kowa{Kowalski K and Rembieli\'nski J 2002 {\it J. Phys. A: Math.
Gen.} {\bf 35} 1405}

\nref\bayen{Bayen F, Flato M, Fronsdal M, Lichnerowicz A and
Sternheimer D 1978 {\it Ann. Phys., NY} {\bf 111}, 61  
\hfill\break
Bayen F, Flato M, Fronsdal M, Lichnerowicz A and
Sternheimer D 1978 {\it Ann. Phys., NY} {\bf 111}, 111} 

\nref\stern1{Sternheimer D 1998 Deformation quantization: twenty years
after {\it Particles, Fields and Gravitation} ed J Rembieli\'nski
(Woodbury, New York: American Institute of Physics) pp 107-145}

\nref\stern2{Dito G and Sternheimer D 2002 Deformation quantization:
genesis, developments and metamorphoses {\it Preprint} math.QA/0201168}

\nref\dito{Dito J 1990 {\it Lett. Math. Phys.} {\bf 20} 125
\hfill\break
Dito J 1993 {\it Lett. Math. Phys.} {\bf 27} 73
\hfill\break
Dito J 2002 Deformation quantization of covariant fields {\it Preprint}
math.QA/0202271}

\nref\anton{Antonsen F 1997 {\it Phys. Rev. D} {\bf 56} 920}

\nref\curt{Curtright T, Fairlie D and Zachos C 1998 {\it Phys. Rev. D}
{\bf 58} 025002
\hfill\break
Curtright T, Uematsu T and Zachos C 2001 {\it J. Math. Phys.} {\bf 42}
2396}

\nref\bor1{Bordemann M and Waldmann S 1998 {\it Commun. Math. Phys.} {\bf 
195} 549}

\nref\bor2{Bordemann M, Neumaier N and Waldmann S 1999 {\it J. Geom.
Phys.} {\bf 29} 199}

\nref\cat{Cattaneo A and Felder G 2000 {\it Commun. Math. Phys.} {\bf 212}
591}

\nref\dqf{Garc\'{\i}a-Compe\'an H, Pleba\'nski J F, Przanowski M and
Turrubiates F J 2000 {\it J. Phys. A: Math. Gen.} {\bf 33} 7935 
\hfill\break
Garc\'{\i}a-Compe\'an H, Pleba\'nski J F, Przanowski M and
Turrubiates F J 2001 {\it Int. J. Mod. Phys. A} {\bf 16} 2533}

\nref\curt2{ Curtright T and Zachos C 2001 {\it Mod. Phys. Lett. A} {\bf 
16} 2381}

\nref\pleb{Pleba\'nski J F 1969 Nawiasy Poissona i komutatory (Toru\'n
{\it Preprint} Nr 69) }

\nref\jacob{Jacobson N 1980 {\it Basic Algebra II} (San Francisco: W. 
H. Freeman and Company)}

\nref\lang{Lang S 1984 {\it Algebra} (Menlo Park, California:
Addison-Wesley Publishing Company, Inc.)}

\nref\fuchs{Fuchs L 1963 {\it Partially Ordered Algebraic Systems}
(London: Pergamon Press)}

\nref\rajw{Rajwade A R 1993 {\it Squares} (Cambridge: Cambridge University
Press)} 

\nref\schar{Scharlau W 1985 {\it Quadratic and Hermitian Forms}
(Berlin: Springer-Verlag)}

\nref\prest{Prestel A and Delzell C N 2001 {\it Positive Polynomials}
(Berlin: Springer Monographs in Mathematics, Springer)}

\nref\artin{Artin E and Schreier O 1926 {\it Abh. Math. Sem. Hamburg}
{\bf 5} 85; {\it The Collected Papers of Emil Artin} 1965 eds S Lang and
J T Tate (USA: Addison-Wesley Publishing Company, Inc.) }

\nref\mie{Bia{\l}ynicki-Birula I, Mielnik B and Pleba\'nski J 1969 
{\it Ann. of Phys.} {\bf 51} 187}

\nref\fed{ Fedosov B 1996 {\it Deformation Quantization and Index Theory} 
(Berlin: Akademie Verlag) }

\nref\kon{Kontsevich M 1997 Deformation quantization of Poisson manifolds
I {\it Preprint} q-alg/9709040}

\nref\prie{Prie\ss-Crampe S 1983 {\it Angeordnete Strukturen: Gruppen,
K\"orper, projektive Ebenen} (Berlin Heidelberg:Springer-Verlag)}

\nref\ruiz{Ruiz J M 1993 {\it The Basic Theory of Power Series} 
(Braunschweig:Vieweg-Verlag)}

\nref\hahn{Hahn H 1907 {\it S.-B. Akad. Wiss. Wien, math.-naturw. Kl.
Abt. IIa} {\bf 116} 601}

\nref\neu{Neumann B H 1949 {\it Trans. Amer. Math. Soc.} {\bf 66} 202}

\nref\mac{MacLane S 1939 {\it Bull. Amer. Math. Soc.} {\bf 45} 888}

\nref\alling{Alling N L 1962 {\it Trans. Amer. Math. Soc.} {\bf 103}
341}

\nref\wilde{De Wilde M and Lecomte P B A 1983 {\it Lett. Math. Phys}
{\bf 7} 487}

\nref\omori{Omori H, Maeda Y and Yoshioka A 1991 {\it Adv. Math.} {\bf 85}
224}

\nref\take{Takesaki M 1979 {\it Theory of Operator Algebras I} (New
York: Springer-Verlag) }

\nref\bord3{ Bordemann M and Waldmann S 1997 {\it Lett. Math. Phys} {\bf
41} 243}



\Title{math.QA/0207181, CINVESTAV-FIS-02/71}
{\vbox{\centerline{}
\vskip -1truecm
\centerline{Uncertainty Relations} 
\centerline{in Deformation Quantization}}}
\vskip -1truecm
\centerline{ M. Przanowski$^{a,b}$\foot{przan@fis.cinvestav.mx} and
F.J. Turrubiates$^a$\foot{fturrub@fis.cinvestav.mx} }
\smallskip
\centerline{\it $^a$Departamento de F\'{\i}sica}
\centerline{\it  Centro de Investigaci\'on y de Estudios Avanzados del
IPN}
\centerline{\it Apdo. Postal 14-740, 07000, M\'exico D.F., M\'exico}
\smallskip
\centerline{\it $^b$Institute of Physics}
\centerline{\it Technical University of \L \'od\'z}
\centerline{\it W\'olcza\'nska 219, 93-005, \L \'od\'z, Poland}
\bigskip
\baselineskip 18pt
\medskip
\vskip 1truecm
\centerline{\it Dedicated to Jerzy Pleba\'nski on the occasion of his
75th birthday}
\vskip 1truecm
\noindent
{\bf Abstract.} Robertson and Hadamard-Robertson theorems on non-negative
definite hermitian forms are generalized to an arbitrary ordered field.
These results are then applied to the case of formal power series fields,
and the Heisenberg-Robertson, Robertson-Schr\"odinger and trace
uncertainty relations in deformation quantization are found.
Some conditions under which the uncertainty relations are minimized are
also given. 

\bigskip



\Date{July, 2002}

\newsec{Introduction} 

The Heisenberg uncertainty relation for canonical observables
$q$ and $p$ is certainly one of the most fundamental results in quantum
mechanics. It was introduced by Heisenberg in 1927 [1] and
mathematically proved by Kennard [2] and Weyl [3]. Later on the
Heisenberg uncertainty relation was generalized to the case of two
arbitrary observables by Robertson [4,5] and Schr\"odinger [6].
In fact in [5,6] an improved version of the Heisenberg 
uncertainty relation has been obtained. Finally, Robertson [7] was able
to extend the previous results to an arbitrary number
of observables. The inequalities found in [7] are called the {\it
Heisenberg-Robertson} and {\it Robertson-Schr\"odinger uncertainty
relations}.

Recently a great deal of interest in uncertainty relations is observed. It
has been shown that they can be used to define {\it squeezed} and
{\it coherent states} and also to generalize this important concepts by
introducing the notion of {\it intelligent states} [8-18].    

It seems to be natural that any theory which would like to describe 
quantum systems should reproduce in some sense the uncertainty relations.
So we expect that it must be also the case in deformation quantization. 

Deformation quantization as introduced by Bayen, Flato, Fronsdal,
Lichnerowicz and Sternheimer [19] and extensively developed during last
years (for review see [20,21]), besides to be a well constructed
mathematical formalism is expected to be an alternative approach to the
description of quantum systems. A big effort in this direction has been
made, e.g. [22-28]. 

The aim of the present paper is to study the uncertainty relations in
deformation quantization. This problem in the case of two observables was
already considered by Curtright and Zachos [29]. 
We are going to extend their results to the the case of an arbitrary
number of observables (real formal power series) and so to obtain in
deformation quantization the Heisenberg-Robertson and
Robertson-Schr\"odinger uncertainty relations and  also the concept of
Robertson-Schr\"odinger intelligent state.   

To deal with uncertainty relations in deformation quantization first one 
should consider the theory of formally real ordered fields and
then apply it to the field of formal power series. This is done in
Sections 2 and 3.    

The importance of the theory of formally real ordered fields in
deformation quantization and especially in the Gel'fand-Naimark-Segal
(GNS) construction was recognized by Bordemann and Waldmann [25]. In our
paper we use extensively the results of their distinguished work.

In Section 4 the proofs of Robertson and Hadamard-Robertson 
theorems for an arbitrary ordered field are given. The results of 
this section are then used in Section 5 to obtain the
Heisenberg-Robertson, Robertson-Schr\"odinger and trace uncertainty
relations. 
Some conditions to minimize the Robertson-Schr\"odinger 
uncertainty relations and to get intelligent states are found in Section
6. These conditions are the deformation quantization analogs of the 
ones introduced by Trifonov [12].

Finally, some conclusion remarks in Section 7 close our paper.

\vskip 1truecm

{\it We would like to dedicate this modest work to our Teacher and Friend,
Professor Jerzy Pleba\'nski who several years ago showed us his
works on Moyal bracket and the beautiful notes from his Polish lectures
entitled ``Nawiasy Poissona i Komutatory" [30]. This was the inspiration
of our interest in deformation quantization}. 

\vskip 1truecm

\newsec{Formally real fields}

In this section we give a brief review of the Artin-Schreier 
theory of formally real fields. For a detailed exposition the reader is
referred to the books by Jacobson [31], Lang [32], Fuchs [33],
Rajwade [34], Scharlau [35], Prestel and Delzell [36]
or to the original paper by Artin and Schreier [37].
 
Let $\IK$ be a field. 

{\bf Definition 2.1}. {\it An ordered field is a pair $(\IK,P)$ where P is
a subset of $\IK$ such that}

{\it (i)} $0 \notin P, \ \   P \cap -P = \emptyset$; 

{\it (ii)} $P + P \subset P, \ \  P \cdot P \subset P$; 

{\it (iii)} $\IK= P \cup \{ 0 \} \cup -P$. 

If $(\IK,P)$ is an ordered field then we say that $\IK$ {\it is ordered
by} $P$ and $P$ is called an {\it order of} $\IK$ or the {\it set of
positive elements of} $\IK$. It is easy to show that if $P$ and $P'$ are
two orders of $\IK$ and $P' \subset P$ then $P' = P$. 

Let $a \not= 0$ be any element of $\IK$. By {\it (iii)} $a \in P$ or $-a
\in P$. Then by {\it (ii)} $a^2=(-a)^2 \in P$. Consequently, if $a_i \in
\IK, i=1,...,n$, then $a_1^2+...+a_n^2=0$ iff $a_i=0$ $\forall$ $i$. 
Now, since $1=1^2 \in P$ one has $1+...+1 \not= 0$ what means that the
characteristic of $\IK$ is $0$.      

One defines the relations $>$ and $\geq$ by: $a>b$ for $a,b \in \IK$ iff
$a-b \in P$; $a \geq b$ iff $a>b$ or $a=b$.
The following properties of the relation $>$ can be easily proved: 
$$
a>0 \ {\rm iff} \  a \in P, 
$$
$$
a>b \ {\rm and} \  b>c \Rightarrow a>c,
$$
$$
a \not= b \Rightarrow a>b  \ {\rm or} \  b>a,
$$
$$
a>b \Rightarrow a+c>b+c \ {\rm for \  any} \ c \in \IK,
$$
$$
a>b \Rightarrow ad>bd \  {\rm for \ any} \ d \in P.
$$

As is used in the real number theory we write $b \leq a$ iff $a \geq b$.

Given an ordered field $(\IK,P)$ the module $| \cdot |$ can be defined by:
$|a| = a$ for $a>0$, $|a|=-a$ for $a \leq 0$. One quickly finds that $|a
b|=|a||b|$ and $|a+b| \leq |a| + |b|$.

{\bf Definition 2.2}. {\it $(\IK,P)$ is called an Archimedean ordered
field if for each $a \in \IK$ there exists an $n \in \IN$ such that
$n 1 >a$}.

An important class of fields called {\it formally real fields} was introduced and analyzed by
Artin and Schreier in their pioneer work [37].
 
{\bf Definition 2.3}. {\it $\IK$ is said to be a formally real field if
$-1$ is not a sum of squares in $\IK$}. 

The classical example of this type of fields is provided by the real
number field $\IR$. Another example fundamental for our further
constructions will be given in the next section.
 
The connection between the ordered fields and the formally real fields is
given by

{\bf Theorem 2.1}. {\it $\IK$ can be ordered iff $\IK$ is a formally
real}. \qed

{\bf Definition 2.4}. {\it A field $\IK$ is called real closed if 

(i) $\IK$ is formally real; 

(ii) any formally real algebraic extension of $\IK$ is equal to $\IK$ }. 

For example the real number field $\IR$ is real closed. 

The following theorems characterize the real closed fields:

{\bf Theorem 2.2}. {\it If $\IK$ is real closed, then $\IK$ has a unique
order $P=(\IK-\{0\})^2:=\{a^2:a \in \IK-\{0\}\}$ }. \qed

{\bf Theorem 2.3}. {\it The following statements are equivalent}: 

{\it (1)} $\IK$ {\it is real closed}. 

{\it (2)} {\it Any polynomial of odd degree with coefficients in $\IK$ has
a root in $\IK$ and there exists an order $P$ of $\IK$ such that any
positive element has square root in $\IK$}.

{\it (3)} {\it $\sqrt{-1} \notin \IK$ and $\IK(\sqrt{-1})$ is
algebraically closed}. \qed

(We use the usual notation in which $\IK(X_1,...,X_n)$ denotes the field
of rational functions in $X_1,...X_n$ with coefficients in $\IK$. So
$\IK(\sqrt{-1})=\IK+\sqrt{-1}\IK$). 

Now the natural question arises if an arbitrary ordered field can be
extended to a real closed one. To answer this question first we give 

{\bf Definition 2.5}. {\it Let $(\IK,P)$ be an ordered field. A field
$\IK'$ is said to be a real closure of $\IK$ relative to $P$ if the
following conditions are satisfied}:

{\it (i)} $\IK'$ {\it is an algebraic extension of $\IK$};

{\it (ii)} $\IK'$ {\it is real closed};

{\it (iii)} $P=(\IK'-\{0\})^2 \cap \IK$ {\it i.e., the unique order
$(\IK'-\{0\})^2$ of $\IK'$ is an extension of $P$}.

Perhaps the most important result in the formally real fields theory is
the theorem due to Artin and Schreier [37] on the existence and
uniqueness of a real closure for any ordered field.

{\bf Theorem 2.4}. {\it Any ordered field $(\IK,P)$ has a real closure
relative to $P$. If $(\IK_1,P_1)$ and $(\IK_2,P_2)$ are ordered
fields and $\IK'_1$ and $\IK'_2$ their respectives closures, then any
isomorphism $f:\IK_1 \rightarrow \IK_2$ such that $f(P_1)=P_2$ can be
uniquely extended to an isomorphism $f':\IK'_1  \rightarrow \IK'_2$ with
$f'((\IK'_1-\{0\})^2)=(\IK'_2-\{0\})^2$ }.  \qed

This theorem will be applied in Section $4$ to prove the generalized
Robertson inequality. 

Finally, we introduce the notion of an exponential valuation of an
arbitrary field $\IK$

{\bf Definition 2.6}. {\it Let $\IK$ be a field. An exponential valuation
of $\IK$ is a mapping $\nu: \IK \rightarrow \IR \cup \{\infty\}$ such that
for all $a,b \in \IK$ 

(i) $\nu(a) = \infty$ $\Leftrightarrow$ $a=0$

(ii) $\nu(ab) = \nu(a) + \nu(b)$

(iii) $\nu(a+b)$ $\geq$ $min\{ \nu(a),\nu(b) \}$ }

Given an exponential valuation $\nu: \IK \rightarrow \IR \cup \{ \infty
\}$ one can define a metric on $\IK$ as follows

\eqn\cero{d_{\nu}(a,b):= \exp \{ -\nu(a-b) \}  \ \ \ \big( \exp \{-\infty 
\}:=0)}

The pair $(\IK,d_{\nu})$ is a metric space, and consequently all
notions known in the theory of metric spaces can be applied in the present
case e.g., a topology ${\cal T}_{\nu}$ defined by the metric $d_{\nu}$,
Cauchy sequences, completeness, etc.

Let $(\IK,P)$ be an ordered field. Then we have a natural topology
${\cal T}_{0}$ on $\IK$: A base of ${\cal T}_{0}$ is the set ${\cal B}$ of
$\varepsilon$-balls, where the {\it $\varepsilon$-ball of the centre in $a
\in \IK$, $B_{\varepsilon}(a)$}, is defined by
$$
B_{\varepsilon}(a):=\{ b \in \IK: |b-a| < \varepsilon \}, \ \ \
0<\varepsilon \in \IK
$$

If the topology ${\cal T}_{\nu}$ on $\IK$ defined by the valuation $\nu:
\IK \rightarrow \IR \cup \{ \infty \}$ is equal to the topology
${\cal T}_{0}$ defined by the order $P$ of $\IK$ then we say that the {\it
valuation $\nu$ is compatible with the ordering of $\IK$}.

An important example when ${\cal T}_{\nu}={\cal T}_{0}$ is considered in
the next section.

\vskip 1truecm
\newsec{Fields of formal power series}

Formal power series play an important role in mathematical physics.
Examples of this are the formal solution of the evolution Schr\"odinger
equation or the Baker-Campbell-Hausdorff formula [38]. But maybe the most
transparent application of formal power series theory can be found in deformation quantization as formulated by 
Bayen et al [19] and developed by Fedosov [39], Kontsevich[40]
and others (see [20,21]). Here the formal power series with respect to the
deformation parameter $\hbar$ arise as the main objects of the
construction.

We give a short exposition of general formal power series field theory. For 
details see [31,33,36,41,42,25,26].

Let $(G,+)$ be an additive abelian group.

{\bf Definition 3.1}. {\it An ordered abelian group is a pair $((G,+),S)$
where $S$ is a subset of $G$ such that} 

{\it (i)} $0 \notin S, \ \   S \cap -S = \emptyset$;

{\it (ii)} $S + S \subset S$;

{\it (iii)} $G= S \cup \{ 0 \} \cup -S$.

We use the symbol $0$ for the neutral element of the group $(G,+)$ as well as for the zero element of a field.

If $g_1, g_2 \in G$ then we say that $g_1<g_2$ ({\it $g_1$ is less than $g_2$}) iff $g_1-g_2 \in S$. So $g \in S$
iff $g<0$ and it means that $S$ consists of elements of $G$ less than the neutral element $0$. We write $g_1 \leq 
g_2$ iff $g_1<g_2$ or $g_1=g_2$.

{\bf Definition 3.2}. {\it Let $((G,+),S)$ be an ordered abelian group
and $\IK$ a field. A formal power series
on $G$ over $\IK$ is a map $a:G \rightarrow \IK$ such that any nonempty
subset of the set $supp\ a := \{g \in G: a(g) \not= 0 \}$ has a least
element.}

Formal power series $a:G \rightarrow \IK$ is usually denoted by $a=\sum_{g
\in G}a_g t^g$ where $a_g:=a(g)$. The
set of all formal power series on $G$ over $\IK$ will be denoted by
$\IK((t^G))$. When $(G,+)$ is the abelian group of
integers $(\IZ,+)$ we simply write $\IK((t))$.  

Addition and multiplication of formal power series $a=\sum_{g \in G}a_g t^g$ 
and $b=\sum_{g \in G}b_g t^g$ are defined as follows

\eqn\uno{a+b =\sum_{g \in G}(a_g + b_g) t^g, \ \ \ \  a b = \sum_{g \in G} \big( \sum_{g_1
\in G}a_{g_1}b_{g-g_1}\big) t^g}

(Note that according to Definition 3.2 both operations are well defined.
In particular for any $g \in G$ the number of non zero elements of the
form $a_{g_1}b_{g-g_1}$, $g_1 \in G$, is finite).

As has been shown by Hahn [43] (and then generalized by Neumann
[44]) {\it the set
$\IK((t^G))$ together with the addition and multiplication defined by Eq.
\uno \ form a field.}

Remind that $(G,+)$ is called a {\it root group} if for any integer $n$ and 
every $g \in G$ there exists $g' \in G$ such that $ng'=g$. We need also
the notion of {\it universal field}.
A field $\IK$ is said to be {\it universal} if every other field $\IK'$ of
the same cardinal number and the same characteristic as $\IK$ is
isomorphic to some subfield of $\IK$.

The following two important theorems have been proved by MacLane [45]
(see also [41]):
 
{\bf Theorem 3.1}. {\it If the coefficient field $\IK$ is algebraically
closed and the ordered abelian
group $G$ is a root group, then the power series field $\IK((t^G))$ is
algebraically closed.} \qed

{\bf Theorem 3.2}. {\it If the coefficient field $\IK$ is algebraically
closed and the ordered abelian group $G$ is
a root group and it contains an element different from the neutral element $0$, then the power series field
$\IK((t^G))$ is universal.} \qed

Suppose that the coefficient field $\IK$ is formally real and is ordered
by $P$. Then $\IK((t^G))$ is a formally real
field and there exists a natural order $P'$ of $\IK((t^G))$ generated by
the order $P$. This order is defined as follows

{\bf Definition 3.3}. {\it If $a=\sum_{g \in G}a_gt^g$, $a_g \in \IK$, and
$g_0$ is the least element of supp a, then $a>0$ iff $a_{g_0}>0$ }.

For the case when the coefficient field $\IK$ is formally real one can
rewrite Theorem 3.1 in the form (see also Alling [46])

{\bf Theorem 3.1$'$}. {\it If the coefficient field $\IK$ is real closed
and the ordered abelian group $G$ is a root group, then the power series
field $\IK((t^G))$ is real closed.} \qed

The fundamental object in the usual deformation quantization construction is 
an associative algebra $(C^{\infty}(M)((\hbar)),*)$ over the complex field
$\IC((\hbar))=\IR((\hbar))+ \sqrt{-1} \IR((\hbar))$. We discuss this
algebra in more details in Section 5. Here we note only that
$C^{\infty}(M)((\hbar))$ denotes the set of formal power series on the
group $\IZ$ with coefficients being smooth complex functions on a
symplectic manifold $M$. (As is used in deformation quantization the
parameter $t$ is denoted by $\hbar$).

However, in the light of Theorems 2.3, 3.1, 3.2 and 3.1$'$ it seems to be
more convenient to deal with the algebra $(C^{\infty}(M)((\hbar^{\IQ})),*)$
over the complex field $\IC((\hbar^{\IQ}))=\IR((\hbar^{\IQ}))+\sqrt{-1}
\IR((\hbar^{\IQ}))$ where $(\IQ,+)$ is the group of rational numbers.     
This conclusion can be also justified from the analytical point of view.

To this end define a valuation $\nu: \IC((\hbar^{\IQ})) \rightarrow \IR
\cup \{\infty\}$ $(or \ \IR((\hbar^{\IQ})) \rightarrow \IR
\cup \{\infty\})$ as follows 
\eqn\unoi{ \nu(a)= min (supp\ a), \ \ \ a \in \IC((\hbar^{\IQ})) \ \ \ (or
\ \IR((\hbar^{\IQ})) )}
Then the metric $d_{\nu}:\IC((\hbar^{\IQ})) \times \IC((\hbar^{\IQ}))
\rightarrow \IR$ (or $\IR((\hbar^{\IQ})) \times \IR((\hbar^{\IQ}))
\rightarrow \IR$) is given by  \cero.

Analogously as it has been done in [25] (Proposition 2) one can prove:

{\bf Proposition 3.1}. {\it $(\IC((\hbar^{\IQ})),d_{\nu})$ and
$(\IR((\hbar^{\IQ})),d_{\nu})$ are complete metric spaces.} \qed

It is also a simple matter to show that the valuation \unoi\ is compatible
with the ordering of $\IR((\hbar^{\IQ}))$ given by Definition 3.3, i.e.,
${\cal T}_{\nu}={\cal T}_{0}$, where the topologies ${\cal T}_{\nu}$ and
${\cal T}_{0}$ are defined in Section 2.

$\{$ {\bf Remark}:  M. Bordemann and S. Waldmann [25] deal with some
subfields of $\IC((\hbar^{\IQ}))$ (or $\IR((\hbar^{\IQ}))$) defined as
follows:

\noindent {\it (1)} The field of {\it formal Newton-Puiseux (NP) series}
$$
\IC \langle\langle \hbar^* \rangle\rangle:= \{ a \in \IC((\hbar^{\IQ})):
\exists \ N \in \IN \ \ N \cdot supp\ a \subset \IZ \}; 
$$
{\it (2)} The field of {\it formal completed Newton-Puiseux (CNP)
series}
$$
\IC \langle\langle \hbar \rangle\rangle:= \{ a \in \IC((\hbar^{\IQ})):
supp\ a \cap [p,q]\ is\ finite\
for\ any \ p,q \in \IQ \},
$$

\noindent and similarly for $\IR \langle\langle \hbar \rangle\rangle$ and
$\IR \langle\langle \hbar^* \rangle\rangle$.
In Proposition 2 of [25] it is shown that
$(\IC\langle\langle \hbar \rangle\rangle,d_{\nu})$ and
$(\IR \langle\langle \hbar \rangle\rangle,d_{\nu})$ are
complete metric spaces. Moreover, $(\IC \langle\langle
\hbar^*\rangle\rangle,d_{\nu})$ (or
$(\IR\langle\langle \hbar^* \rangle\rangle,d_{\nu})$) is dense in
$(\IC \langle\langle \hbar \rangle\rangle,d_{\nu})$ (or
$(\IR \langle\langle \hbar \rangle\rangle,d_{\nu})$). Then in Theorem 1 of
[25] it is proved that
both fields, $\IC \langle\langle \hbar^*\rangle\rangle$ and
$\IC \langle\langle \hbar \rangle\rangle$, are algebraically closed   
($\IR \langle\langle \hbar^* \rangle\rangle$ and $\IR
\langle\langle \hbar \rangle\rangle$ are real closed).

It is evident that $(\IC \langle\langle \hbar \rangle\rangle,d_{\nu})$
and, consequently,
$(\IC \langle\langle \hbar^* \rangle\rangle,d_{\nu})$ are not dense metric
spaces in
$(\IC((\hbar^{\IQ})),d_{\nu})$. However, since
$(\IC((\hbar^{\IQ})),d_{\nu})$ is complete and the
field $\IC((\hbar^{\IQ}))$ is algebraically closed then
$\IC((\hbar^{\IQ}))$ can be applied to the GNS construction in deformation
quantization analogously as it is in the case of $\IC
\langle\langle \hbar \rangle\rangle$ [25,26]. $\}$

\vskip 1truecm
\newsec{Robertson and Hadamard-Robertson theorems for an arbitrary
ordered field}

The well known Heisenberg uncertainty relation between two canonical 
observables admits several generalizations. One of them was given by
Robertson [5] and Schr\"odinger [6].
These results were then generalized to an arbitrary number of
observables by Robertson [7]. Recently a revival interest of the
important Robertson work can be observed ([12-17] and the references given
therein).

In this section we are going to generalize Robertson's results to an 
arbitrary formal real ordered field. Let $(\IK,P)$ be a formally real
ordered field and $\IK^c:=\IK(i)=\IK+i \IK$, $i \equiv \sqrt{-1}$, its
complexification. 

Let $V$ be a vector space over $\IK^c$. 

{\bf Definition 4.1}. {\it A hermitian form on $V$ is a map $\phi:V
\times V \rightarrow \IK^c$ satisfying the following properties.

(i) $\phi(c_1v_1+c_2v_2,w)=\overline{c_1}\phi(v_1,w)+\overline{c_2}\phi(v_2,w)$

(ii) $\phi(v,c_1w_1+c_2w_2)=c_1\phi(v,w_1)+c_2\phi(v,w_2)$

(iii) $\overline{\phi(v,w)}= \phi(w,v)$

$\forall \  v_1,v_2,w_1,w_2,v,w \in V$, $\forall \  c_1,c_2 \in \IK^c$}

In this paper the overbar denotes the complex conjugation.

({\bf Note:} A map $\psi:V \times V \rightarrow \IK^c$ is said 
to be a {\it sesquilinear form} if it satisfies {\it (i)} and {\it (ii)}
[32]). 

Hermitian form $\phi:V \times V \rightarrow \IK^c$ is said to be {\it
positive definite} if $\phi(v,v)>0 \ for \ all \ nonzero \ v \in V$; and
it is said to be {\it non-negative definite} if $\phi(v,v) \geq
0$ $\ \forall \  v \in V$.   

Suppose that $dim V=n$. Denote by $e_1=(1,0,...,0),...,e_n=(0,...,0,1)$ 
the natural basis of $V$. Let $v=\sum_{j=1}^n v_je_j$  be any vector of
$V$. Then from Definition 4.1 one gets
\eqn\dos{\phi(v,v)=\sum_{j,k=1}^n\phi_{jk}\overline{v_j}v_k, \ \ \
\overline{\phi_{jk}}=\phi_{kj}}
where $\phi_{jk}:=\phi(e_j,e_k)$. 

We can write $\phi_{jk}=a_{jk}+ i b_{jk}$, $a_{jk},b_{jk} \in \IK$. From
\dos \ if follows that
$a_{jk}=a_{kj}$ and $b_{jk}=-b_{kj}$. So the $n \times n$ matrix
($\phi_{jk}$) over $\IK^c$ is {\it hermitian}, the matrix $(a_{jk})$ over
$\IK$ is {\it symmetric} and the matrix $(b_{jk})$ over $\IK$ is {\it
skew-symmetric}.

Now we are in a position to prove a generalization of the Robertson theorem to an arbitrary formally real ordered
field.

{\bf Theorem 4.1} ({\it Robertson}). {\it With the notation as 
above, let $\phi:V \times V \rightarrow \IK^c$ be a non-negative definite
hermitian form on $V$. Then
$det(a_{jk}) \geq det(b_{jk})$. If $\phi$ is positive definite then
$det(a_{jk}) > det(b_{jk})$. If $det(a_{jk})=0$ then $det(b_{jk})=0$}.

{\bf Proof:}

Let $v=\sum_{j=1}^n v_je_j$ be any vector in $V$. Write $K^c \ni v_j=x_j +
i y_j$, $x_j,y_j \in \IK$.
Then $\phi(v,v)= \sum_{j,k=1}^n a_{jk}(x_jx_k+y_jy_k)-2 \sum_{j,k=1}^n
b_{jk} x_jy_k$. Letting $y_j=0$, one
quickly finds that 

$\phi(v,v) \geq 0, \forall \  v \in V \Rightarrow \sum_{j,k=1}^n
a_{jk}x_jx_k \geq 0 \ \ \  \forall \  x_j \in \IK$.
 
Consequently, $\sum_{j,k=1}^n a_{jk}x_jx_k$ is a non-negative definite
quadratic form and in particular it follows that $det(a_{jk}) \geq 0$.

The $n \times n$ matrix $(b_{jk})$ is skew-symmetric. Hence if $n$ is an odd number then $det(b_{jk})=0$ and the
theorem holds. Thus we assume that $n$ is an even number $n=2m$.

The proof is divided into two parts.

(1) $det(a_{jk})>0$.

Here we follow Robertson [7] (see also [16]).

The matrices $(a_{jk})$ and $(c_{jk}):=i (b_{jk})$ are hermitian. One
can find a $2m \times 2m$ matrix $D$ over
$\IK$, with $det(D) \not= 0$, such that the transformed matrix
$(a'_{jk}):=D^T(a_{jk})D$, where $D^T$ denotes the
transposed matrix of $D$, is a diagonal matrix with all its diagonal
elements positive. According to Theorem 2.4
we can extend the ordered field $(\IK,P)$ to its closure $(\IK',P')$. So
without any loss of generality we assume from
the very beginning that $(\IK,P)$  {\bf is real closed}.
With this assumption and by Theorem 2.3 the matrix $D$ over $\IK$ can be
found such that the matrix $(a'_{jk})$ is the unit matrix $1$.
It is obvious that the matrix $(c'_{jk}):=D^T (c_{jk})D$ is still
hermitian. Therefore, there exists a $2m \times 2m$ unitary matrix $U$
over $\IK^c$ ($U^{\dagger}U=1$, $U^{\dagger}:=\overline{U}^T$) such that
the transformed matrix $(c''_{jk}):=(DU)^{\dagger}(c_{jk})(DU)$ is
diagonal. Moreover, the transformed matrix
$(a''_{jk}):=U^{\dagger}(a'_{jk})U = 1$. Hence, finally we get
\eqn\tres{(c''_{jk})=(DU)^{\dagger}(c_{jk})(DU)=diag(\lambda_1,...,\lambda_{2m}),
\ \ \ (a''_{jk})=(DU)^{\dagger}(a_{jk})(DU)=1}     

where $\lambda_1,...,\lambda_{2m} \in \IK$ are the solutions of the
characteristic equation

\eqn\cuatro{det((c'_{jk}) - \lambda 1) = 0}

Since the matrix $(c'_{jk})$ is skew-symmetric then if $\lambda$ is a
solution of the characteristic equation \cuatro\ then $-\lambda$ is also a
solution of this equation.

Therefore the matrix $(c''_{jk})$ is of the form

\eqn\cinco{(c''_{jk})=diag(\lambda_1,-\lambda_1,...,\lambda_m,-\lambda_m),
\ \ \ \lambda_1,...,\lambda_m \in \IK}

By Eqs. \tres \ and \cinco \ one quickly finds that the transformed matrix
$(\phi''_{jk})$ of the hermitian form $\phi$ reads

\eqn\seis{(\phi''_{jk}):=(DU)^{\dagger}(\phi_{jk})(DU)=(a''_{jk}) +
(c''_{jk})=diag(1+\lambda_1, 1-\lambda_1,...,1+\lambda_m,1-\lambda_m)}

Since $\phi$ is a non-negative definite hermitian form then

\eqn\siete{1 \pm \lambda_k \geq 0, \ \ \  k=1,...,m}

From \tres \ and \cinco \ one gets

\eqn\ocho{det(a_{jk})=(detD)^{-2}, \ \ \ \ \ 
det(c_{jk})=(detD)^{-2} (-1)^m \lambda_1^2...\lambda_m^2}

But $det(c_{jk})=det(ib_{jk})=i^{2m}det(b_{jk})=(-1)^mdet(b_{jk})$.

Substituting this relation into \ocho, employing also the fact that by
\siete \  $\lambda_1^2...\lambda_m^2 \leq 1$ we obtain that the inequality
$det(a_{jk}) \geq det(b_{jk})$ holds true.

\noindent Observe that if $\phi$ is a positive definite hermitian form
then $det(\phi_{jk})>0$ $\Rightarrow$ $det(a_{jk})>0$. Moreover, in \siete
\ one has the strict inequalities $1 \pm \lambda_k > 0$ and consequently,
we obtain the strict inequality $det(a_{jk})>det(b_{jk})$.

This completes the first part of the proof. Consider now the second part
when: 

(2) $det(a_{jk})=0$

Then it follows that also $det(\phi_{jk})=0$.

There exists an unitary matrix $U$ over $\IK^c$ such that

\eqn\nueve{(\phi'_{jk}):=U^{\dagger}(\phi_{jk})U=
diag(\phi_1,...,\phi_q,0,...,0), \ \ \  q<2m,  \ \ \ \
\phi_1,...,\phi_q>0}

Define now

\eqn\diez{(\phi'_{jk}(x)):=diag(\phi_1,...,\phi_q,x,...,x),  \ \ \ \ \  x
\geq 0}

It is evident that the hermitian form $\phi(x)$ given by the matrix

\eqn\once{(\phi_{jk}(x))=U(\phi'_{jk}(x))U^{\dagger}, \ \ \ \ \  x \geq
0}

is positive definite for every $x>0$. Moreover,
$(\phi_{jk}(0))=(\phi_{jk})$ i.e., $\phi(0)=\phi$.

We split $(\phi_{jk}(x))$ as before

\eqn\doce{\phi_{jk}(x)=a_{jk}(x)+ib_{jk}(x), \ \ \ \ \  x \geq 0,
\ \ \  a_{jk}(x), b_{jk}(x) \in \IK}
$$
a_{jk}(x)=a_{kj}(x), \ b_{jk}(x)=-b_{kj}(x); \ \ \ a_{jk}(0)=a_{jk}, \
b_{jk}(0)=b_{jk}
$$

Since $det(\phi_{jk}(x))>0$ \ \ $\forall \ x>0$ then also $det(a_{jk})>0$
\ \ $\forall \ x>0$ and by the first part (1) of the proof one has

\eqn\trece{det(a_{jk}(x)) > det(b_{jk}(x)) \ \ \ \ \ \forall \ x > 0}

From \nueve, \  \diez, \  \once, \  \doce \  and the fact that
$det(a_{jk})=0$ it follows that

\eqn\catorce{det(a_{jk}(x))=det(a_{jk})+ \sum_{l=1}^{r \leq 2m} d_lx^l =
\sum_{l=1}^{r \leq 2m}d_lx^l}
$$
det(b_{jk}(x)) = det(b_{jk}) + \sum_{l=1}^{s \leq 2m} f_lx^l,   \ \ \ \ \
d_l,f_l \in \IK 
$$ 

Consequently, by \trece \ and \catorce

\eqn\quince{\sum_{l=1}^{p \leq 2m}g_lx^l - det(b_{jk}) > 0,  \ \ \ \forall
\  x>0}
$$
\sum_{l=1}^{p \leq 2m}g_lx^l:=\sum_{l=1}^{r \leq 2m} d_lx^l -
\sum_{l=1}^{s \leq 2m} f_lx^l.
$$

Since $(b_{jk})$ is a skew-symmetric matrix over the formally real field
$\IK$ then $det(b_{jk}) \geq 0$. Hence, $g_l \not= 0$ for some $l$. 

We will show that $det(b_{jk}) = 0$. 

Suppose that $det(b_{jk})>0$. The inequality \quince \ yields

\eqn\dieciseis{ \bigg( \sum_{l=1}^{p \leq 2m} |g_l|x^l - det(b_{jk})
\bigg) > 0, \ \ \ \ \ \forall \ x > 0}

Without any loss of generality one can assume that all $g_l \not= 0$. Put
then

\eqn\diecisiete{ x = min \bigg( {det(b_{jk}) \over 2p|g_1|},
\sqrt{det(b_{jk}) \over 2p|g_2|},..., \root p\of{det(b_{jk})
\over 2p|g_p|} \bigg) }

Remeber that, as it has been pointed out in the first part (1) of our proof,
without any loss of generality one can consider $\IK$ to be a real
closed field. So \diecisiete \ is well defined by the Theorem 2.3.

Substituting $x$ given by \diecisiete \ into \dieciseis \ we infer that

\eqn\dieciocho{ \bigg({det(b_{jk}) \over 2} - det(b_{jk}) \bigg) > 0
\Rightarrow det(b_{jk}) < 0}

This contradicts the assumption: $det(b_{jk})>0$. Consequently,
$det(b_{jk}) = 0$ and the proof is complete.  \qed

$\{${\bf Remark}: Note that using analytical methods a different proof of
the second part (2) of Theorem 4.1 can be given. Namely, taking the limit
of both sides of the inequality \quince\ when $x \rightarrow 0^+$ one
immediately gets

\eqn\diecinueve{\lim_{x\to 0^+} \bigg( \sum_{l=1}^{p \leq 2m}g_lx^l -
det(b_{jk}) \bigg) \geq 0.}

As $\lim_{x\to 0^+} \big( \sum_{l=1}^{p \leq 2m}g_lx^l \big)=0$ and
$det(b_{jk}) \geq 0$ we obtain $det(b_{jk}) = 0$. $\}$  

From the proof of Theorem 4.1 (especially see (4.5) and (4.7)) we find 
that for $n=2$ i.e. $m=1$, the following corollary holds:

{\bf Corollary 4.1}. {\it If $n=2$ then $det(a_{jk})=det(b_{jk})$ iff
$det(\phi_{jk})=0$.} \ \  \qed

One can prove a useful lemma which will be employed to
generalize the Hadamard-Robertson theorem. 

Keeping the notation as above one has
  
{\bf Lemma 4.1}. {\it Let $\phi:V \times V \rightarrow$ $\IK^c$ be a
non-negative definite hermitian form on $V$. Then $det(a_{jk}) \geq
det(\phi_{jk})$. 
Equality $det(a_{jk})=det(\phi_{jk})$ holds iff $det(a_{jk})=0$ or
$(\phi_{jk})= (a_{jk})$.}

{\bf Proof:} 

As before $\phi_{jk}=a_{jk}+ib_{jk}$, where $a_{jk}=a_{kj}$ and
$b_{jk}=-b_{kj}$ are elements of $\IK$. We put $c_{jk}:=ib_{jk}$. 

If $det(a_{jk})=0$ then $det(\phi_{jk})=0$ and the lemma is valid.

Let $det(a_{jk})>0$. Suppose $dimV=n$. Analogously as in Theorem 4.1 one
can choose a $n \times n$ matrix $D$ over $\IK$ such that $(a'_{jk}):= D^T
(a_{jk})D=1$. Obviously the matrix $(c'_{jk}):=D^T(c_{jk})D$ is hermitian
and skew-symmetric. Then an unitary $n \times n$ matrix can be found for
which 
$$
(c''_{jk}):=U^{\dagger}(c'_{jk})U=diag(\lambda_1,...,\lambda_n) \  and \
(a''_{jk}):=U^{\dagger}(a'_{jk})U=1
$$
where $\lambda_1,...,\lambda_n \in \IK$ are the solution of the
characteristic equation \cuatro.
The $n \times n$ matrix $(c'_{jk})$ is skew-symmetric and as before it
follows that if $\lambda$ is a solution of the characteristic equation
then $- \lambda$ is also a solution. 

Hence in the case of even $n$, $n=2m$, we have
$(c''_{jk})=diag(\lambda_1, -\lambda_1,..., \lambda_m,-\lambda_m)$, and in
the case when $n$ is odd, $n=2m+1$, the matrix $(c''_{jk})=diag(\lambda_1,
-\lambda_1,..., \lambda_m,-\lambda_m,0)$. Consequently,
$(\phi''_{jk})=diag(1+\lambda_1,1-\lambda_1,...,1+\lambda_m,1-\lambda_m)$
for $n=2m$, and 
$(\phi''_{jk})=diag(1+\lambda_1,1-\lambda_1,...,1+\lambda_m,1-\lambda_m,1)$
for $n=2m+1$.

Since $\phi$ is non-negative definite then $1 \pm \lambda_k \geq 0$ for
all $k$. Therefore,
$det(\phi''_{jk})=(1-{\lambda_1}^2)...(1-{\lambda_m}^2) \leq 1$. So
$det(\phi''_{jk}) \leq det(a''_{jk})$, and the equality $det(\phi''_{jk})
= det(a''_{jk})$ holds iff $\lambda_1=...=\lambda_m=0$ i.e., iff
$(c''_{jk})=0$. This yields $det(\phi_{jk}) \leq det(a_{jk})$ and the
equality $det(\phi_{jk}) = det(a_{jk})$ holds iff $(c_{jk})=0.$ The proof
is complete.  \qed

To obtain a generalization of the Heisenberg uncertainty principle to any
formally real field it is necessary to generalize first the
Hadamard-Robertson theorem [7].

{\bf Theorem 4.2} ({\it Hadamard-Robertson}). {\it  
Let $\phi: V \times V \rightarrow \IK^c$ be a non-negative definite
hermitian form on a vector space $V$of dimension $n$ over $\IK^c$. Then,  

(i) $\phi_{11}...\phi_{nn} \geq det(a_{jk}) \geq det(\phi_{jk})$,
$\ \ \ \phi_{11}...\phi_{nn} \geq det(a_{jk}) \geq det(b_{jk})$

(ii) $\phi_{11}...\phi_{nn}= det(a_{jk}) = det(\phi_{jk)}$
$\Leftrightarrow$ $\phi_{kk}=0$ for some $k$, or $(\phi_{jk})=(a_{jk})$ is
diagonal. 

(iii) $\phi_{11}...\phi_{nn} = det(b_{jk})$  $\Leftrightarrow$ 
$\phi_{kk}=0$ for some $k$ or $(a_{jk})$ is diagonal and
$det(b_{jk})=det(a_{jk})$.}

{\bf Proof:}
 
{\it (i)}  From Theorem 4.1 and Lemma 4.1 one has:
$det(a_{jk}) \geq det(b_{jk})$ and $det(a_{jk}) \geq det(\phi_{jk})$,
respectively. Hence it remains only to prove that
$\phi_{11}...\phi_{nn} \geq det(a_{jk})$.
But as $\phi_{kk}=a_{kk}$ for $k=1,...,n$ this inequality is equivalent to

\eqn\veinte{a_{11}...a_{nn} \geq det(a_{jk})}

From the assumption that the hermitian form $\phi:V \times V \rightarrow
\IK^c$ is non-negative definite it follows that the quadratic form

\eqn\veintiuno{\sum_{j,k=1}^na_{jk}x_jx_k, \ \ \ \ \  x_j \in \IK}

\noindent is also non-negative definite. Keeping this in mind we prove
\veinte\ by induction with respect to the dimension of $V$. For $dim V=1$
the inequality \veinte\ holds trivially.
Assume that \veinte\ is valid for $dim V=n-1$, $n \geq 2$. Let now $dim
V=n$. We can find a $n \times n$ orthogonal matrix $R$ over $\IK$ of the
form

\eqn\veintidos{R=\pmatrix{r_{11}&...&r_{1,n-1}&0\cr .&...&.&.\cr
r_{n-1,1}&...&r_{n-1,n-1}&0\cr 0&...&0&1 \cr}, \ \ \ R^TR=1}

such that

\eqn\veintitres{(a'_{jk}):=R^T(a_{jk})R =
\pmatrix{\lambda_1&0&...&0&a'_{1n} \cr 0&\lambda_2 & ....&0&a'_{2n} \cr.&.
&...&.&.\cr 0&0&...&\lambda_{n-1}&a'_{n-1,n} \cr a'_{1 n}&a'_{2
n}&....&a'_{n-1,n}&a'_{nn} \cr},  \ \ \ \lambda_1,...,\lambda_{n-1} \geq
0, \ \ \ a'_{nn}=a_{nn}}

One quickly finds that

\eqn\veinticuatro{det(a'_{jk})=det(a_{jk})=\lambda_1...\lambda_{n-1}a_{nn}
- (a'_{11})^2 \lambda_2...\lambda_{n-1}- \lambda_1(a'_{2n})^2
\lambda_3...\lambda_{n-1} -}
$$
... - \lambda_1
\lambda_2...\lambda_{n-2}(a'_{n-1,n})^2  \leq det(A_{n-1})a_{nn} 
$$

where $A_{n-1}$ is the $(n-1) \times (n-1)$ matrix over $\IK$ defined by
$$
A_{n-1} := \pmatrix{a_{11}&...&a_{1,n-1}  \cr a_{21}&...&a_{2,n-1} \cr
.&...&.\cr a_{n-1,n}& ...& a_{n-1,n-1} \cr},
\ \ \ detA_{n-1}=\lambda_1...\lambda_{n-1} 
$$
Since the quadratic form \veintiuno\ is non-negative definite, then the
quadratic form $\sum_{j,k=1}^{n-1}a_{jk}x_jx_k$, $x_j \in \IK,$ is also
non-negative definite. Consequently, the inductive assumption gives
$$
a_{11}...a_{n-1,n-1} \geq det(A_{n-1})
$$

Substituting this into \veinticuatro\ one gets \veinte\ and the proof of
{\it (i)} is complete.

{\it (ii)}  $\Leftarrow$  \ \ \ If $\phi_{kk}=0$ for some $k$ or
$(\phi_{jk})=(a_{jk})$ is diagonal then 

\eqn\veinticinco{\phi_{11}...\phi_{nn} = det(a_{jk}) =
det(\phi_{jk})}

$\Rightarrow$ \ \ \ Assume that \veinticinco\ holds. Hence, from Lemma 4.1
we conclude that
$$
det(\phi_{jk})=det(a_{jk})=0 \ \ \ or  \ \ \ (\phi_{jk})=(a_{jk})
$$
Obviously, $det(\phi_{jk})=0$ with \veinticinco\ imply that $\phi_{kk}=0$
for some $k$. Suppose then that \veinticinco\ is valid and
$det(\phi_{jk}) > 0$. Now $(a_{jk})=(\phi_{jk})$ and in \veintitres\
$\lambda_1,...,\lambda_{n-1},a_{nn} >0$. So from \veinticuatro\ it follows
that the equality $det(a_{jk})=det(A_{n-1}) \cdot a_{nn}$ holds iff,
$a'_{1n}=...=a'_{n-1,n}=0$.

This last condition by \veintidos\ and \veintitres, is equivalent to
$a_{1n}=...=a_{n-1,n}=0$.

Analogous considerations for $A_{n-1}$,....etc.; lead to the conclusion
that \veinticinco\ with $det(\phi_{jk})>0$ imply
$(a_{jk})=(\phi_{jk})=diag(\phi_{11},...,\phi_{nn})$, $\phi_{kk}>0$ for
$k=1,...,n$.

{\it (iii)}  The proof is straightforward keeping in mind that
$det(b_{jk}) \geq 0$ and employing {\it (i)} and {\it (ii)}.   \qed

Finally, we would like to generalize to an arbitrary formally real
ordered field $(\IK,P)$ an interesting uncertainty relations for the trace
of the matrix $(\phi_{jk})$ (Trifonov [17]).

{\bf Proposition 4.1}. {\it For any non-negative definite hermitian form
$\phi:V
\times V \rightarrow \IK^c$, the following inequality holds 
\eqn\veintiseis{Tr(\phi_{jk}) \geq {2 \over n-1} \sum_{j<k}^{n}|b_{jk}|}
for every $n$, where $n=dim V$. If n is even, n=2m, then also 

\eqn\veintiocho{Tr(\phi_{jk}) \geq 2 \sum_{j=1}^m |b_{j,m+j}|} 
}

{\bf Proof:} Assume $j \not= k$. We start with the obvious relation
$$
(a_{jj}+a_{kk})^2 \geq 4a_{jj}a_{kk}
$$
From the Hadamard-Robertson Theorem 4.2 we have 
$$
a_{jj} a_{kk} \geq b_{jk}^2
$$
Consequently,
\eqn\veintinueve{a_{jj}+a_{kk} \geq 2|b_{jk}|}
$$ 
a_{jj}+a_{kk}=2|b_{jk}| \Leftrightarrow a_{jk}=0  \ {\rm and} \ 
a_{jj}=a_{kk}=|b_{jk}|
$$
Using the relation
$$
Tr(\phi_{jk})=Tr(a_{jk})={1 \over n-1} \sum_{j<k}^n(a_{jj}+a_{kk})
$$ 
and \veintinueve\  one gets that \veintiseis\ holds true.

If $n=2m$ we can write
$$
Tr(\phi_{jk})=Tr(a_{jk})=\sum_{j=1}^m(a_{jj} + a_{m+j,m+j})
$$ 
This with \veintinueve\ give \veintiocho\ and the proposition is proved.
\qed


\newsec{Uncertainty relations in deformation quantization}

Deformation quantization was introduced as an alternative approach to the
description of quantum systems. In the fundamental work by Bayen,
Flato, Fronsdal, Lichnerowicz and Sternheimer [19] it is
suggested that quantization should be understood ``... as a deformation
of the structure of the algebra of classical observables, rather than a
radical change in the nature of the observables". 
This construction is realized by a deformation of the usual product
algebra of smooth functions on the phase space and then by a deformation
of the Poisson algebra. 

To be more precise: Let $(M,\omega)$ be a symplectic manifold ($\omega$ 
denotes the symplectic form on $M$), and let $C^{\infty}(M)((\hbar))$ be a
vector space over $\IC((\hbar))$ of the formal power series
\eqn\treinta{f=\sum_{k=-N}^{\infty}f_k(x)\hbar^k}
where $f_k(x)$ are complex smooth functions on $M$, $f_k \in
C^{\infty}(M)$.

{\bf Definition 5.1}. [19,39,25] {\it Deformation quantization on
$(M,\omega)$ is an associative algebra $(C^{\infty}(M)((\hbar)),*)$ over
the field $\IC((\hbar))$, where the associative product $*$, called
star-product, is given by}  
\eqn\treintaiuno{f * g = \sum_{k=0}^{\infty}C_k(f,g)\hbar^k, \ \ \ f,g \in
C^{\infty}(M)((\hbar))}
{\it with $C_k$, $k \geq 0$, being bidifferential operators such that
$C_k(C^{\infty}(M) \times C^{\infty}(M)) \subset C^{\infty}(M) \ \forall
k$, $C_k(1,f)=C_k(f,1)=0$ for $k \geq 1$, $C_0(f,g)=fg$,
$C_1(f,g)-C_1(g,f)=i\hbar \{f,g\}$ and $\{
\cdot,\cdot \}$ stands for the Poisson bracket}.  

It has been proved [47,48,39] that deformation quantization exists on each
symplectic manifold. Even more, recently Kontsevich [40] proved the
existence of star product for an arbitrary Poisson manifold.
Perhaps the most transparent construction of star product on an
arbitrary symplectic manifold has been given by Fedosov [39] in
terms of the geometry of the formal Weyl algebra bundles. For our purpose
is not necessary to consider Fedosov's construction in more detail. 

As it has been pointed out in Section 3 it seems to be natural to extend
the associative algebra $(C^{\infty}(M)((\hbar)),*)$ over the field
$\IC((\hbar))$ to $(C^{\infty}(M)((\hbar^{\IQ})),*)$ over the field
$\IC((\hbar^{\IQ}))$. In what follows we deal with such an extended
deformation quantization.

To proceed further we need the definition of positive functionals and
states in deformation quantization. These concepts are fundamental in the
GNS construction developed by M. Bordemann et al [25,26]
and so seem to be basic to relate deformation quantization with quantum
mechanics.   

Analogously as in the theory of $C^*$-algebras one has [25,49]:

{\bf Definition 5.2}. {\it A $\IC^{\infty}((\hbar^{\IQ}))$ linear
functional $\rho:C^{\infty}(M)((\hbar^{\IQ}))
\rightarrow \IC((\hbar^{\IQ}))$ is said to be positive if
$$
\rho(\overline{f}*f) \geq 0 \ \ \ \forall f \in
C^{\infty}(M)((\hbar^{\IQ}))
$$
A positive linear functional $\rho$ is called a state if $\rho(1)=1$}.

One can easily check that if a linear functional $\rho$ is positive then
\eqn\treintaitres{\overline{\rho(f*g)} = \rho(\overline{g}*\overline{f})}
and the {\it Cauchy-Schwarz inequality}
\eqn\treintaitresb{\rho(\overline{f}*g)\overline{\rho(\overline{f}*g)}
\leq \rho(\overline{f}*f)\rho(\overline{g}*g)}
holds true.
In particular taking in \treintaitres\  $g=1$ we get
\eqn\treintaicuatro{\overline{\rho(f)} =\rho(\overline{f}).}
Consequently, if $\overline{f}=f$ then $\rho(f) \in \IR((\hbar^{\IQ}))$. 

From \treintaitres\ and \treintaicuatro\ it follows that
$$
\rho(\overline{f*g} - \overline{g} * \overline{f})=0.
$$ 
This condition is satisfied for any positive functional iff
\eqn\treintaicuatrob{\overline{f*g}=\overline{g} * \overline{f} \ \ \
\forall f,g \in C^{\infty}(M)((\hbar^{\IQ}))} 
Note that it is always possible to construct a star product which
satisfies \treintaicuatrob\ [39,50].

Another fundamental concept in the GNS construction and employed in the
present paper to describe intelligent states (Section 6) is that of the
{\it Gel'fand ideal}

{\bf Definition 5.3}. {\it Let $\rho:C^{\infty}(M)((\hbar^{\IQ}))
\rightarrow \IC((\hbar^{\IQ}))$ be a positive linear functional. Then the 
subspace ${\cal J}_{\rho}$ of $C^{\infty}(M)((\hbar^{\IQ}))$
$$
{\cal J}_{\rho}:= \{f \in C^{\infty}(M)((\hbar^{\IQ})):
\rho(\overline{f}*f)=0 \}
$$
is called the Gel'fand ideal of $\rho$.}

It can be easily shown that by \treintaitres\ and \treintaitresb\ ${\cal 
J}_{\rho}$ is a left ideal of $C^{\infty}(M)((\hbar^{\IQ}))$, i.e. if $f
\in {\cal J}_{\rho}$ then $g*f \in {\cal J}_{\rho}$ $\forall g \in
C^{\infty}(M)((\hbar^{\IQ}))$ and

\eqn\treintaicuatroc{\rho(\overline{f}*g)=0=\rho(g*f) \ \ \ \forall g \in
C^{\infty}(M)((\hbar^{\IQ})) } 

Let $\rho:C^{\infty}(M)((\hbar^{\IQ})) \rightarrow \IC((\hbar^{\IQ}))$
be a positive linear functional. Define the sesquilinear form
$\phi:C^{\infty}(M)((\hbar^{\IQ})) \times C^{\infty}(M)((\hbar^{\IQ}))
\rightarrow \IC((\hbar^{\IQ}))$ by 
\eqn\treintaicinco{\phi(f,g):= \rho(\overline{f}*g)  \ \ \  f,g \in
C^{\infty}(M)((\hbar^{\IQ})). }

(For the definition of a sesquilinear form see the note after Definition
4.1)

From \treintaitres\ one quickly finds that
\eqn\treintaiseis{\overline{\phi(f,g)}=\phi(g,f).}

It means that $\phi$ is a Hermitian form on
$C^{\infty}(M)((\hbar^{\IQ}))$. Moreover, since
$\phi(f,f)=\rho(\overline{f}*f) \geq 0$ $\forall f \in
C^{\infty}(M)((\hbar^{\IQ}))$ then $\phi$ defined by \treintaicinco\ is a
{\bf non-negative definite hermitian form}.  

Now we are at the position to obtain uncertainty relations in deformation
quantization. To this end, let $X_1,...,X_n \in
C^{\infty}(M)((\hbar^{\IQ}))$ satisfy the reality conditions
$\overline{X_j}=X_j, \ j=1,...,n$ (i.e., $X_j$ are {\it observables}) and
let $\rho:C^{\infty}(M)((\hbar^{\IQ})) \rightarrow \IC((\hbar^{\IQ}))$ be
a state. Define {\it deviations from the mean} as follows

\eqn\treintaisiete{\delta X_j:= X_j - \rho(X_j)}

Since $\overline{X_j}=X_j$  and $\rho$ is a state then by \treintaicuatro\
one gets 

\eqn\treintaiocho{\overline{\delta X_j}=\delta X_j.}

It is also evident that $\rho(\delta X_j)=0$.
Take
$$
f:= \sum_{j=1}^n v_j \delta X_j, \ \ \ \ v_j \in \IC((\hbar^{\IQ}))
$$
Then from \treintaicinco\ and \treintaiocho\ we have 
$$
\phi(f,f)=\rho \big( \sum_{j=1}^n \overline{v_j} \delta X_j * \sum_{k=1}^n
v_k \delta X_k \big) = \sum_{j,k=1}^n \overline{v_j}v_k \rho(\delta X_j *
\delta X_k)=\sum_{j,k=1}^n \overline{v_j}v_k \phi (\delta X_j, \delta X_k) 
$$
Define
\eqn\treintainueve{\phi_{jk}:= \rho(\delta X_j * \delta X_k)=\phi(\delta 
X_j , \delta X_k), \ \ \ \phi_{jk} \in \IC((\hbar^{\IQ})). }
From \treintaiseis\ it follows that $\overline{\phi_{jk}}=\phi_{kj}$. 
Since $\phi(f,f) \geq 0$ then 
\eqn\cuarenta{ \sum_{j,k=1}^n  \phi_{jk}\overline{v_j}v_k \geq 0, \ \ \
\forall v_j \in \IC((\hbar^{\IQ})).}
Consequently, the $n \times n$ hermitian matrix $(\phi_{jk})$ over
$\IC((\hbar^{\IQ}))$ determines a non-negative hermitian form \cuarenta.  

We can use now the results of Section 4. 

First, as before we write $\phi_{jk}=a_{jk}+ib_{jk}$, $a_{jk}, b_{jk} \in
\IR((\hbar^{\IQ}))$. From \treintainueve\ and \treintaisiete\ one gets 

\eqn\cuarentaiuno{a_{jk}= {1 \over 2} \rho \big((\delta X_j * \delta X_k +
\delta X_k * \delta X_j) \big)= {1 \over 2} \rho\big((X_j * X_k + X_k *
X_j)\big) - \rho(X_j) \rho(X_k)=a_{kj} } 
$$
b_{jk}= -{i \over 2} \rho \big((\delta X_j * \delta X_k -
\delta X_k * \delta X_j) \big)= {\hbar \over 2} \rho ( \{X_j,X_k \}_*
)=-b_{kj}
$$
where $\{ X_j,X_k \}_*:= {1 \over i \hbar}(X_j * X_k - X_k * X_j)$. 
In analogy to quantum mechanics and statistics the $n \times n$
symmetric matrix $(a_{jk})$ over $\IR((\hbar^{\IQ}))$ can be called the
{\it dispersion} or {\it covariance matrix}. A diagonal element
$a_{jj}=\rho(X_j * X_j)-(\rho(X_j))^2$  which we denote also by $(\Delta
X_j)^2$ is the {\it variance of $X_j$}, and $\Delta X_j = \sqrt{a_{jj}}$
is the {\it uncertainty in $X_j$} (or {\it standard deviation of $X_j$}). 
The element $a_{jk}$ for $j \not= k$ is the {\it covariance of
$X_j$ and $X_k$}. 

Having all that the Theorem 4.1 leads to the following {\it
Robertson-Schr\"odinger uncertainty relation} in deformation quantization:

\eqn\cuarentaidos{det \bigg( {1 \over 2} \rho( \delta X_j * \delta X_k +
\delta X_k * \delta X_j) \bigg) \geq det \bigg( {\hbar \over 2} \rho( \{X_j,X_k 
\}_* ) \bigg).}

In particular for two observables $X_1$ and $X_2$ we get

\eqn\cuarentaitres{\Delta X_1 \Delta X_2 \geq {1 \over 2} \sqrt{ (\hbar
\{ X_1,X_2 \}_*)^2 + \big( \rho(X_1*X_2+X_2*X_1)-2\rho(X_1)\rho(X_2)
\big)^2} }

This is the deformation quantization analogue of the well known in quantum
mechanics uncertainty relation given by Robertson [5] and Schr\"odinger
[6].
The relation \cuarentaitres\ has been found recently by Curtright and
Zachos [29]. However, their result seems to be derived in the spirit
of a strict deformation quantization which makes use of Wigner
function and not for the formal deformation quantization in the sense
of Bayen et al [19] considered in the present paper.    

Another uncertainty relation in deformation quantization which we call the
{\it Heisenberg-Robertson uncertainty relation} follows immediately from
the Hadamard-Robertson theorem (Theorem 4.2), and it reads

\eqn\cuarentaicuatro{(\Delta X_1)^2...(\Delta X_n)^2 \geq det \bigg(
{\hbar \over 2} \rho( \{ X_j,X_k \}_*  ) \bigg) }

Finally, employing the Proposition 4.1 one gets the {\it trace uncertainty 
relation}

\eqn\cuarentaicinco{(\Delta X_1)^2+...+(\Delta X_n)^2 \geq {\hbar \over
n-1} \sum_{j<k}^n |\rho(\{X_j,X_k \}_*)| }


\newsec{Intelligent states in deformation quantization}

In quantum mechanics the states that minimize the Heisenberg-Robertson or
the Robertson-Schr\"odinger uncertainty relations play an important
role in the theory of coherent and squeezed states and they are called
{\it Heisenberg-Robertson} or {\it Robertson-Schr\"odinger intelligent
states}, ({\it minimum uncertainty states, correlated coherent states})
[8-13,16]. It seems to be reasonable to extend these notions to
deformation quantization. Thus we have 

{\bf Definition 6.1}. {\it A state $\rho:C^{\infty}(M)((\hbar^{\IQ}))
\rightarrow \IC((\hbar^{\IQ}))$ is said to be a Heisenberg-Robertson
intelligent state for $X_1,...,X_n$ if
\eqn\cuarentaiseis{(\Delta X_1)^2...(\Delta X_n)^2 = det \bigg(
{\hbar \over 2} \rho( \{ X_j,X_k \}_*  ) \bigg).}
If
\eqn\cuarentaisiete{det \bigg( {1 \over 2} \rho( \delta X_j * \delta X_k +
\delta X_k * \delta X_j) \bigg) = det \bigg( {\hbar \over 2} \rho(
\{X_j,X_k \}_* ) \bigg)}
then $\rho$ is called a Robertson-Schr\"odinger intelligent state for
$X_1,...,X_n$.} 

\noindent From Theorems $4.1$ and $4.2$ one can easily obtain that
$$
\cuarentaiseis\ \Rightarrow \cuarentaisiete.
$$   
Hence every Heisenberg-Robertson intelligent state is also a
Robertson-Schr\"odinger intelligent state.

To have a deeper insight into the Robertson-Schr\"odinger intelligent
states we prove some conditions under which \cuarentaisiete\ is satisfied.  

Our results are the deformation quantization versions of the
propositions found by Trifonov in the case of quantum
mechanics (Propositions: 1 and 3 of [12]).  

Observe that by Theorem 4.1 if $det \bigg( {1 \over 2} \rho( \delta X_j * 
\delta X_k + \delta X_k * \delta X_j) \bigg) = 0$ then also $det \bigg(
{\hbar \over 2} \rho(\{X_j,X_k \}_* ) \bigg) =0$. Hence, $det \bigg( {1
\over 2} \rho( \delta X_j * \delta X_k + \delta X_k
* \delta X_j) \bigg) = 0$ is a sufficient condition for $\rho$ to be a
Robertson-Schr\"odinger intelligent state for $X_1,...,X_n$. In the case
when the number $n$ of observables $X_j$ is odd this condition is also
necessary.

We can prove

{\bf Proposition 6.1}. {\it Let $\rho:C^{\infty}(M)((\hbar^{\IQ}))
\rightarrow \IC((\hbar^{\IQ}))$ be a state and $a_{jk}:={1 \over
2}\rho(\delta X_j * \delta X_k + \delta X_k * \delta X_j)$, $j,k=1,..,n$.
Then $det(a_{jk})=0$ iff there exist $x_1,...,x_n \in \IR((\hbar^{\IQ}))$
such that $\sum_{j=1}^n|x_j|>0$ and
\eqn\cuarentaiocho{\rho\bigg( \sum_{j=1}^n x_j \delta X_j * \sum_{k=1}^n
x_k \delta X_k \bigg) =0}
i.e., $\sum_{j=1}^n x_j \delta X_j$ is an element of the Gel'fand ideal
${\cal J}_{\rho}$ of $\rho$. }

{\bf Proof:} (Compare with [12]).

Assume that $det(a_{jk})=0$.
Then there exists a $n \times n$ orthogonal matrix $R=(r_{jk})$ over
$\IR((\hbar^{\IQ}))$, $R^TR=1$, such that 
$$
R^T(a_{jk})R=R^T \bigg( {1 \over 2} \rho(\delta X_j * \delta X_k + \delta
X_k * \delta X_j) \bigg) R=diag(\lambda_1,...,\lambda_{q-1},0...,0), \ \ \
2 \leq q  \leq n.
$$
Hence 
$$
\rho \bigg( \sum_{j=1}^n r_{jq} \delta X_j * \sum_{k=1}^nr_{kq}
\delta X_k \bigg)=0.
$$

Denoting $x_j:=r_{jq} \in \IR((\hbar^{\IQ}))$ one gets \cuarentaiocho .
This completes the first part of the proof.

Assume now that there exist $x_1,...,x_n \in \IR((\hbar^{\IQ}))$ such that
$\sum_{j=1}^n|x_j|>0$ and \cuarentaiocho\ holds. Choose a $n \times n$
matrix $D=(d_{jk})$ over $\IR((\hbar^{\IQ}))$ such that $d_{j1}=x_j$,
$j=1,...,n$, and $det D \not= 0$.

Consider the transformed matrix $(a'_{jk})=D^T(a_{jk})D$.

We have
$$
a'_{1l}={1 \over 2} \rho \bigg( \sum_{j=1}^n x_j \delta X_j * \sum_{k=1}^n
d_{kl} \delta X_k + \sum_{k=1}^n d_{kl} \delta X_k * \sum_{j=1}^n
x_{j} \delta X_j \bigg), \ \ \ l=1,...,n.
$$
Since $\sum_{j=1}^nx_j \delta X_j \in {\cal J}_{\rho}$ then by
\treintaicuatroc\ 
$$
\rho \bigg( \sum_{j=1}^n x_j \delta X_j * g \bigg) = 0 = \rho \bigg( g *
\sum_{j=1}^n x_j \delta X_j \bigg) \ \ \ \forall g \in
C^{\infty}(M)((\hbar^{\IQ})).  
$$

Therefore, $a'_{1l}=0$ for $l=1,...,n$ and consequently, $det(a'_{jk})=0$.
But $det(a'_{jk})=(det D)^2 det(a_{jk})$ with $det D \not=0.$ This yields
$det(a_{jk})=0$.
The proof is complete.  \qed
    
To find another sufficient condition that a given state $\rho$ be a
Robertson-Schr\"odinger intelligent state for $X_1,...,X_n$ we deal with
the case when $n$ is an even number, $n=2m$. Thus we have $X_1,...,X_{2m} 
\in C^{\infty}(M)((\hbar^{\IQ}))$ such that $\overline{X_j}=X_j$,
$j=1,...,2m$. Let $\delta X_j$ be deviations from the mean as in
\treintaisiete\. Introduce the following objects 
\eqn\cuarentainueve{\delta A_{\alpha}:={1 \over 2} (\delta X_{\alpha} + 
i \delta X_{\alpha+m}) }
$$
\overline{\delta A_{\alpha}}={1 \over 2} (\delta X_{\alpha} -
i \delta X_{\alpha+m}), \ \ \ \alpha=1,...,m
$$

With all that one has

{\bf Proposition 6.2}. {\it If there exists a linear transformation 
\eqn\cincuenta{\delta A'_{\alpha} = \sum_{\beta=1}^m  \bigg( u_{\alpha
\beta} \delta A_{\beta} + v_{\alpha \beta} \overline{\delta A_{\beta}}
\bigg)} 
$$
\overline{\delta A'_{\alpha}} = \sum_{\beta=1}^m  \bigg(
\overline{v_{\alpha \beta}} \delta A_{\beta} + \overline{u_{\alpha \beta}}
\overline{\delta A_{\beta}} \bigg); \ \ \  
u_{\alpha \beta}, v_{\alpha \beta} \in \IC((\hbar^{\IQ})),  \  \alpha,
\beta=1,...,m
$$
such that 
\eqn\cincuentaiuno{ det \pmatrix{(u_{\alpha \beta})&(v_{\alpha \beta} )
\cr  &  \cr(\overline{v_{\alpha \beta}})&(\overline{u_{\alpha \beta}})\cr}
\not= 0.}
and 
\eqn\cincuentaidos{\rho(\overline{\delta A'_{\alpha}}* \delta
A'_{\alpha})=0; \ \ \ \alpha=1,...,m}
($\delta A'_{\alpha}$ belongs to the Gel'fand ideal ${\cal
J}_{\rho}$), then \cuarentaisiete\ is satisfied i.e., $\rho$ is a
Robertson-Schr\"odinger intelligent state for $X_1,...,X_{2m}$}.

{\bf Proof:} (Compare with [12]).

Following \cuarentainueve\ define 
$$
\delta X'_{\alpha}:=(\delta A'_{\alpha} + \overline{\delta A'_{\alpha}})
$$
$$
\delta X'_{\alpha+m}:= -i (\delta A'_{\alpha} - \overline{\delta
A'_{\alpha}}),  \ \ \ \alpha =1,...,m. 
$$
Obviously $\overline{\delta X'_j}=\delta X'_j$, $j=1,...,2m$ and one can
easily check that 

\eqn\cincuentaitres{\delta X'_j = \sum_{k=1}^{2m} d_{jk} \delta X_k}
where under \cuarentainueve , \cincuenta\ and \cincuentaiuno\ the $2m
\times 
2m$ matrix $(d_{jk})$ over $\IR((\hbar^{\IQ}))$is non-singular,
$det(d_{jk}) \not=0$.
Straightforward calculations under the assumption \cincuentaidos\ lead to
the relation
$$
det \bigg( {1 \over 2} \rho( \delta X'_j * \delta X'_k +
\delta X'_k * \delta X'_j) \bigg) = det \bigg( {\hbar \over 2} \rho(
\{X'_j,X'_k \}_* ) \bigg).
$$
Consequently, by \cincuentaitres\ the equation \cuarentaisiete\ holds
true. \qed

Employing Corollary 4.1 for the case of two observables one can easily
prove the next proposition

{\bf Proposition 6.3} {\it A state $\rho$ is a Robertson-Schr\"odinger
intelligent state for $X_1$, $X_2$ iff there exist $u_1,u_2 \in
\IC((\hbar^{\IQ}))$ such that $u_1 \delta X_1 + u_2 \delta X_2 \in {\cal 
J}_{\rho}$. } \qed

Robertson-Schr\"odinger intelligent states for two observables in terms of
Moyal star product and Wigner functions have been considered in [24,29].


\newsec{Concluding remarks}

In this paper we have obtained uncertainty relations in deformation
quantization formalism. To achieve this, first it was necessary to study a
general theory of formal real ordered fields and to apply it to the
case of formal power series. Having done all that we were able to
generalize the Robertson and Hadamard-Robertson theorems to be valid 
for an arbitrary ordered field. This allowed us to formulate several
uncertainty relations and to introduce the concept of intelligent states
in deformation quantization.
Of course further investigations in this direction are needed. In
particular one should consider some concrete set of observables and get
examples of the corresponding intelligent states. 

It is expected that the results of the present paper will give a better
understanding of the relations between quantum mechanics and deformation
quantization.

\vskip 1.5truecm

\centerline{\bf Acknowledgements}

This paper was supported by the CONACYT Project: 32427-E. 


\listrefs

\end